\documentclass[11pt,reqno]{amsart}

\usepackage[margin=1.1in]{geometry}
\usepackage{amsmath,amssymb,amsthm,mathtools}
\usepackage{enumitem}
\usepackage[dvipsnames]{xcolor}
\usepackage{tikz}
\usepackage{float}
\usepackage[colorlinks=true,linkcolor=blue!60!black,citecolor=red!60!black,urlcolor=blue!60!black,pagebackref=true]{hyperref}
\numberwithin{equation}{section}
\newtheorem{theorem}{Theorem}[section]
\newtheorem{lemma}[theorem]{Lemma}
\newtheorem{proposition}[theorem]{Proposition}

\theoremstyle{definition}
\newtheorem{definition}[theorem]{Definition}
\theoremstyle{remark}
\newtheorem{remark}[theorem]{Remark}

\newcommand{\R}{\mathbb R}
\newcommand{\dd}{\,\mathrm d}
\newcommand{\MT}{\mathrm{MT}}

\newcommand{\weakto}{\rightharpoonup}

\begin{document}


\title[Logarithmic Trudinger--Moser suprema in dimension $N$]{A uniform lower bound for logarithmic Trudinger--Moser suprema in dimension $N\ge3$}
 \author{Jos\'e Francisco de Oliveira}
 \author[J.M. do \'O]{Jo\~ao Marcos do \'O}
 \author[R.C. Ponciano]{Raon\'{\i} Cabral Ponciano}

\address{Department of Mathematics, Federal University of Piau\'i, 64049-550 Teresina, PI, Brazil}
 \email{\href{mailto:jfoliveira@ufpi.edu.br}{jfoliveira@ufpi.edu.br}}
\address[J.M. do \'O]{Department of Mathematics,
	Federal University of Para\'{\i}ba
	\newline\indent 
	58051-900, Jo\~ao Pessoa-PB, Brazil}
\email{\href{mailto:jmbo@mat.ufpb.br}{jmbo@mat.ufpb.br}}
\address[R.C. Ponciano]{Department of Mathematics,
	Federal University of ABC
	\newline\indent
	09280-560, Santo Andr\'e-SP, Brazil}
\email{\href{mailto:raoni.ponciano@ufabc.edu.br}{raoni.ponciano@ufabc.edu.br}}

\subjclass[2020]{46E35, 26D10, 35J20, 35B33}
\keywords{Trudinger--Moser inequality; logarithmic weights; extremal functions; concentration level; test functions}
\begin{abstract}
Let $B$ be the unit ball of $\R^N$, $N\ge2$, and let $\MT_N(\alpha_{N,\beta},\beta)$ denote the critical Trudinger--Moser supremum of Calanchi and Ruf \cite{CalanchiRuf15,CalanchiRufN15} for radial functions in the Sobolev space with logarithmic weight $(-\log|x|)^{\beta(N-1)}$, $0\le\beta<1$. We prove that $\MT_N(\alpha_{N,\beta},\beta)>(2N-\frac25)|B|$ for every $\beta\in[0,1)$ and $N\ge 3$. This lower bound is uniform in $\beta$ and lies strictly above the upper bound $|B|\big(1+e^{1+\frac12+\cdots+\frac1{N-1}}\big)$ for the concentration level obtained by Carleson--Chang \cite{CarlesonChang86} for $\beta=0$ and, for $\beta>0$, by Roy \cite{Roy16} for $N=2$ and by Roy \cite{RoyN19} and Nguyen \cite{Nguyen19} for $N\ge 3$. As a consequence, we prove that $\MT_N(\alpha_{N,\beta},\beta)$ is attained for every $\beta\in[0,1)$ and $N\ge3$, which was previously known only for $\beta>0$ small enough. For $N=2$, the existence of extremals for every $\beta\in [0, 1)$ was recently obtained in \cite{Oliveira26}.
\end{abstract}
\maketitle

\section{Introduction and main results}\label{sec:introduction}

Let $\Omega\subset\R^N$, $N\ge2$, be a smooth bounded domain and let $W^{1,N}_0(\Omega)$ be the Sobolev space given by the completion of $C^\infty_0(\Omega)$ with respect to the Dirichlet norm $\|\nabla u\|_{N}=\big(\int_\Omega|\nabla u|^N\dd x\big)^{1/N}$. The classical Trudinger--Moser inequality \cite{Moser71,Pohozaev65,Trudinger67,Yudovich61} states that
\begin{equation}\label{eq:Moser}
\MT_{N  
}(\alpha)=\sup_{\substack{u\in W^{1,N}_0(\Omega),\,\|\nabla u\|_N\le1}}\int_{\Omega} e^{\alpha|u|^{\frac{N}{N-1}}}\dd x\;
\begin{cases}
<\infty, & \text{if }\alpha\le\alpha_N,\\
=\infty, & \text{if }\alpha>\alpha_N,
\end{cases}
\end{equation}
where $\alpha_N=N\omega_{N-1}^{\frac1{N-1}}$ and $\omega_{N-1}$ is the measure of the unit sphere in $\R^N$. Inequality \eqref{eq:Moser} gives the sharp exponential integrability of functions in $W^{1,N}_0(\Omega)$, the borderline case of the Sobolev embedding, and has become a fundamental tool in geometric analysis and in the study of nonlinear differential equations with critical exponential growth; see, for instance, \cite{BairdFardounRegbaoui06,ChangYang87,ChangYang88,deFigueiredodoORuf11,Moser71} and the references therein. It has also been extended and refined in many directions, including higher-order derivatives, improved inequalities of Adimurthi--Druet type, the whole space $\R^N$, weighted Sobolev spaces and fractional dimensions; see, without any claim of completeness, \cite{Adams88,AdimurthiDruet04,CalanchiRuf15,deOliveiradoO14,LamLuZhang19} and the references therein. 

An intriguing question concerning the supremum \eqref{eq:Moser} is whether $\MT_{N}(\alpha_N)$ admits an extremal function in the critical case $\alpha=\alpha_N$. This question has been answered affirmatively through the contributions of several authors. The starting point was the work of Carleson and Chang \cite{CarlesonChang86}, who proved the existence of an extremal function when $\Omega=B$ is a ball in $\R^N$, $N\ge2$. Subsequently, Struwe \cite{Struwe88} obtained the existence of a maximizer when $\Omega$ is close to a ball in $\R^2$, and Flucher \cite{Flucher92} extended this result to arbitrary smooth bounded domains in $\R^2$. Finally, Lin \cite{Lin96} established the existence for arbitrary smooth bounded domains in $\R^N$.

In this paper, we are mainly interested in a Trudinger--Moser--type inequality with logarithmic weights in the Sobolev norm due to Calanchi--Ruf \cite{CalanchiRuf15,CalanchiRufN15}. 

For $0\le\beta<1$, let 
$$w_\beta(x)=\big(-\log |x|\big)^{\beta(N-1)}, \quad x\in B\setminus\{0\},
$$ where $B$ is the unit ball centered at the origin in $\R^N$ and consider the weighted Dirichlet norm
\[
\|u\|_{w_\beta}=\Big(\int_B|\nabla u|^Nw_\beta\dd x\Big)^{\frac1N},\;\;\mbox{for}\;\; u\in C^\infty_{0,\mathrm{rad}}(B).
\]
Let $W^{1,N}_{0,\mathrm{rad}}(B,w_\beta)$ denote the completion of $C^\infty_{0,\mathrm{rad}}(B)$ with respect to $\|\cdot\|_{w_\beta}$. Calanchi and Ruf \cite{CalanchiRuf15,CalanchiRufN15} proved that
\begin{equation}\label{eq:sup}
\MT_N(\alpha,\beta)=\sup_{\substack{u\in W^{1,N}_{0,\mathrm{rad}}(B,w_\beta),\,\|u\|_{w_\beta}\le1}}\int_B e^{\alpha|u|^{\gamma_{N,\beta}}}\dd x\;
\begin{cases}
<\infty, & \text{if }\alpha\le\alpha_{N,\beta},\\
=\infty, & \text{if }\alpha>\alpha_{N,\beta},
\end{cases}
\end{equation}
where $\gamma_{N,\beta}=\frac{N}{(N-1)(1-\beta)}$ and $\alpha_{N,\beta}=N\big[(1-\beta)\,\omega_{N-1}^{\frac1{N-1}}\big]^{\frac1{1-\beta}}$. For $\beta=0$, by the P\'olya--Szeg\"{o} inequality, $\MT_N(\alpha,0)$ coincides with the supremum $\MT_{N}(\alpha)$ in \eqref{eq:Moser} with $\Omega=B$, for every $\alpha>0$. Extensions of these inequalities to unbounded domains, to fractional dimensions and to the $k$-Hessian operator can be found in \cite{Aouaoui25,DoOOliveiraPonciano26,XueZhangZhu25}. Nguyen \cite[Theorem~1.1]{Nguyen19} gave a simpler proof of the Calanchi--Ruf inequality \eqref{eq:sup} based on the monotonicity
\begin{equation}\label{eq:mono}
\MT_N(\alpha_{N,\beta'},\beta')\le\MT_N(\alpha_{N,\beta},\beta),\;\;\text{for}\quad0\le\beta\le\beta'<1.
\end{equation}
The existence of maximizers for the critical suprema $\MT_N(\alpha_{N,\beta},\beta)$ in \eqref{eq:sup} has proven to be an interesting challenge. Following a concentration-compactness approach, as in the case $\beta=0$ in \cite{CarlesonChang86,Flucher92,Lin96}, the effect of the weight $w_{\beta}$ on the concentration level of the associated critical functional has already been estimated in \cite{RoyN19} and \cite{Nguyen19}. However, in order to rule out the concentration of a maximizing sequence, a lower bound for the supremum $\MT_N(\alpha_{N,\beta},\beta)$ strictly above this concentration level is required. Following \cite{CarlesonChang86,Roy16} (see also \cite{Nguyen19}), we say that a sequence $(u_k)\subset W^{1,N}_{0,\mathrm{rad}}(B,w_\beta)$ is a normalized concentrating sequence if
\begin{equation}\label{concentraion-sequence}
 \|u_k\|_{w_\beta}=1,\;\; u_k\weakto0\;\;\mbox{weakly in}\;\; W^{1,N}_{0,\mathrm{rad}}(B,w_\beta)\;\;\mbox{and}\;\;\int_{B\setminus B_r}|\nabla u_k|^Nw_\beta\dd x\to0\;\;\mbox{for all}\;\;r\in(0,1), 
\end{equation}
 where $B_r$ is the ball of radius $r$ centered at the origin. Let us denote  by
 \[
 J_{N,\beta}(u)=\int_Be^{\alpha_{N,\beta}|u|^{\gamma_{N,\beta}}}\dd x, \;\;\mbox{for}\;\; u\in W^{1,N}_{0,\mathrm{rad}}(B,w_\beta)
  \]
 the corresponding critical functional. Then we define the concentration level
\[
\Sigma_{N,\beta}=\sup\left\{\limsup_{k\to\infty}J_{N,\beta}(u_k):(u_k)\text{ is a normalized concentrating sequence}\right\}.
\]
To obtain the existence of an extremal function in the case $\beta=0$, the key idea of Carleson and Chang \cite{CarlesonChang86} was to estimate the concentration level
\begin{equation}\label{eq:CC-level}
\Sigma_{N,0}\le |B|\big(1+e^{1+\frac12+\cdots+\frac1{N-1}}\big)
\end{equation}
and then to show, by means of an explicit test function, that
\begin{equation}\label{Sup-Value}
\MT_N(\alpha_N,0)>|B|\big(1+e^{1+\frac12+\cdots+\frac1{N-1}}\big).
\end{equation}
More precisely, arguing by contradiction, they proved that if $\MT_N(\alpha_N,0)$ is not attained, then every normalized maximizing sequence is, up to a subsequence, a normalized concentrating sequence. By \eqref{eq:CC-level}, this gives $\MT_N(\alpha_N,0)\le\Sigma_{N,0}$, which contradicts \eqref{Sup-Value}. Subsequently, this method, together with  Lions' concentration-compactness principle \cite{Lions85},  has been applied to a wide range of extremal problems; see, for example, \cite{CsatoRoy15,deOliveiradoO15,deOliveiradoORuf20,DoOOliveiraPonciano26,DongLamLu18,JiangXuZhangZhu25}.

Roy \cite[Theorem~3.3]{RoyN19} extended \eqref{eq:CC-level} to all $\beta\in [0,1)$, that is,
 \begin{equation}\label{eq:CC-level-geral}
\Sigma_{N,\beta}\le |B|\big(1+e^{1+\frac12+\cdots+\frac1{N-1}}\big),\;\; \mbox{for all}\;\; \beta\in [0, 1).
\end{equation}
Independently, Nguyen \cite[Theorem~1.2]{Nguyen19} proved that  $\beta\mapsto\Sigma_{N,\beta}$ is nonincreasing on $[0,1)$, which, together with \eqref{eq:CC-level}, implies
\begin{equation}\label{eq:concentration}
  \Sigma_{N, \beta'}\le\Sigma_{N,\beta}\le \Sigma_{N, 0}\le |B|\big(1+e^{1+\frac12+\cdots+\frac1{N-1}}\big),\;\;\text{for}\quad0\le\beta\le\beta'<1.
\end{equation}
However, the strict inequality
\begin{equation}\label{Sup-Value-geral}
\MT_N(\alpha_{N,\beta},\beta)>|B|\big(1+e^{1+\frac12+\cdots+\frac1{N-1}}\big)
\end{equation}
has not been fully established so far for $N\geq 3$ and the whole range $\beta\in [0, 1)$. In fact, by test function computations, Roy  \cite{RoyN19} proved it for $N\ge3$ only for $\beta>0$ small enough and concluded that $\MT_N(\alpha_{N,\beta},\beta)$ is attained in this range, while Nguyen \cite{Nguyen19} obtained the same result by employing the monotonicity \eqref{eq:mono}, together with the continuity of $\beta\mapsto\MT_N(\alpha_{N,\beta},\beta)$ at $\beta=0$ \cite[Lemma~3.1]{Nguyen19}. In dimension $N=2$, this problem was first studied by Roy \cite{Roy16}; we refer to \cite[Section~1]{Oliveira26} for a detailed account of the two-dimensional case. Recently, de Oliveira \cite{Oliveira26} proved the strict inequality for the full range $\beta\in [0,1)$ when $N=2$, and concluded the existence of extremals. In this paper, we extend the approach of \cite{Oliveira26} to every dimension $N\ge3$. Our main result is the following uniform lower bound.

\begin{theorem}\label{thm:main}
Let $N\ge3$. Then, for every $0\le\beta<1$,
\begin{equation}\label{eq:main}
\MT_N(\alpha_{N,\beta},\beta)>\big(2N-\frac25\big)|B|>|B|\big(1+e^{1+\frac12+\cdots+\frac1{N-1}}\big).
\end{equation}
\end{theorem}
As in \cite{Oliveira26}, the strict inequality \eqref{eq:main} excludes the concentration of maximizing sequences, and the concentration-compactness alternative \cite[Lemma~2.1]{RoyN19} gives the following.

\begin{theorem}\label{thm:attainment}
For every $N\ge3$ and every $0\le\beta<1$, the supremum $\MT_N(\alpha_{N,\beta},\beta)$ is attained.
\end{theorem}
Theorem~\ref{thm:attainment} removes the smallness assumption on $\beta$ in \cite[Theorem~1.3]{RoyN19} and \cite[Theorem~1.3]{Nguyen19}. Together with \cite{Oliveira26}, it shows that $\MT_N(\alpha_{N,\beta},\beta)$ is attained for every $N\ge2$ and every $\beta\in[0,1)$.

Let us comment on the relation between our approach and that of \cite{Oliveira26}. As in \cite{Oliveira26}, the monotonicity \eqref{eq:mono} allows us to bound $\MT_N(\alpha_{N,\beta},\beta)$ from below, uniformly in $\beta\in[0,1)$, by an explicit integral obtained by letting $\beta\to1^-$ in the normalized test functions (see Proposition~\ref{prop:uniform}); no Trudinger--Moser inequality for $\beta=1$ is involved in this step. Section~\ref{sec:profiles} carries out this reduction in every dimension $N\ge2$, for a class of profiles which contains the polygonal profiles of \cite{Oliveira26}. The main new ingredient is the choice of the test profile. In \cite{Oliveira26}, where $N=2$, an explicit three-node polygonal profile is used, and the resulting integral is estimated with constants adapted to that profile. Here we need, for every $N\ge3$, a profile whose limiting integral can be estimated by a single argument, valid for all $N\ge3$, and compared with the concentration bound $1+e^{1+\frac12+\cdots+\frac1{N-1}}$, which grows linearly in $N$. 
\begin{remark}\label{rem:classical}
For $\beta=0$, Theorem~\ref{thm:main} gives explicit test functions for the strict inequality of Carleson and Chang \cite{CarlesonChang86} between the classical supremum on the unit ball and the concentration level \eqref{eq:CC-level}. Since our bounds are uniform in $\beta$, they are weaker than what can be obtained for $\beta=0$ alone. In the opposite direction, Kim \cite[Theorem~3.5]{Kim00} proved an upper bound for the classical supremum with explicit constants in every dimension; for $N=3$, with the constants given in \cite[Lemma~3.7]{Kim00}, it yields
\[
\MT_3(\alpha_3,0)\le3^{\frac32}e^{\frac38+\frac{16\sqrt5}{27e}}\,|B|<12.31\,|B|.
\]
By \eqref{eq:mono}, this is also an upper bound for $\MT_3(\alpha_{3,\beta},\beta)$, so that Theorem~\ref{thm:main} gives $\frac{28}{5}\,|B|=5.6\,|B|<\MT_3(\alpha_{3,\beta},\beta)<12.31\,|B|$ for every $\beta\in[0,1)$.
\end{remark}

The paper is organized as follows. In Section~\ref{sec:profiles}, we introduce a class of test functions and obtain a lower bound for $\MT_N(\alpha_{N,\beta},\beta)$ that is uniform in $\beta\in[0,1)$. In Section~\ref{sec:test-profile}, we choose a test function modeled on those of Carleson and Chang \cite{CarlesonChang86}, for which this lower bound can be computed almost explicitly. In Section~\ref{sec:proofs}, we prove Theorems~\ref{thm:main} and~\ref{thm:attainment}.

\section{\texorpdfstring{Test functions in dimension $N$}{Test functions in dimension N}}\label{sec:profiles}
Throughout this section $N\ge2$ is fixed, and $\beta$ always denotes a fixed number in $[0,1)$.

\begin{definition}\label{def:profile}
Let $0<t_1<T$. A profile with nodes $t_1$ and $T$ is a Lipschitz continuous, nondecreasing function $\phi:[0,\infty)\to[0,\infty)$ such that $\phi(t)=(y_1/t_1)\,t$ on $[0,t_1]$ for some $y_1>0$, and $\phi(t)=\phi(T)$ for $t\ge T$.
\end{definition}

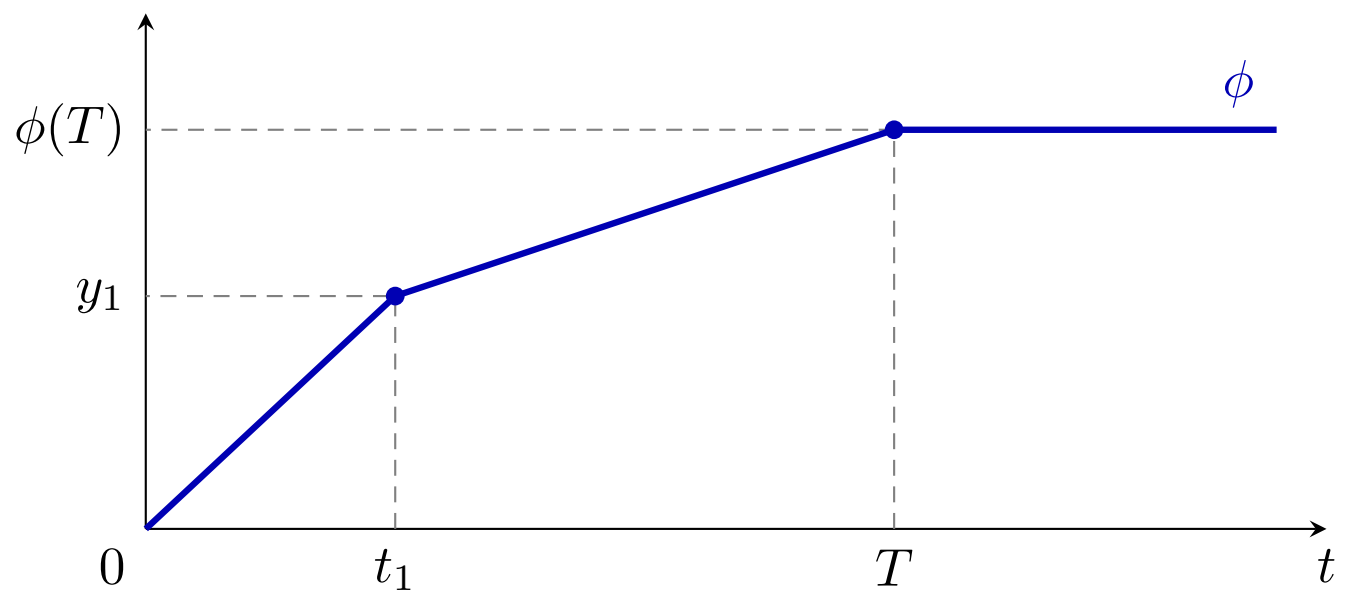
\begin{figure}[ht]
  \centering
  \begin{tikzpicture}[scale=1.1, >=stealth, font=\small]
    \def\tone{1.5} \def\yone{1.4}
    \def\T{4.5}    \def\yT{2.4}
    \def\xmax{6.8}
    \colorlet{profile}{blue!70!black}
    \draw[->] (0,0) -- (\xmax+0.3,0) node[below] {$t$};
    \draw[->] (0,0) -- (0,3.1);
    \draw[densely dashed, gray] (\tone,0) -- (\tone,\yone) -- (0,\yone);
    \draw[densely dashed, gray] (\T,0) -- (\T,\yT) -- (0,\yT);
    \draw[very thick, profile] (0,0) -- (\tone,\yone) -- (\T,\yT) -- (\xmax,\yT)
      node[above left] {$\phi$};
    \fill[profile] (\tone,\yone) circle (1.6pt) (\T,\yT) circle (1.6pt);
    \node[below]      at (\tone,0) {$t_1$};
    \node[below]      at (\T,0)    {$T$};
    \node[left]       at (0,\yone) {$y_1$};
    \node[left]       at (0,\yT)   {$\phi(T)$};
    \node[below left] at (0,0)     {$0$};
  \end{tikzpicture}
  \caption{A profile with nodes $t_1$ and $T$.}
  \label{fig:profile}
\end{figure}

For instance, the polygonal profiles used in \cite{Oliveira26} for $N=2$ are profiles in this sense.

\begin{lemma}\label{lem:test}
Let $\phi$ be a profile with nodes $t_1$ and $T$, and set
\begin{equation}\label{eq:test}
v_{\beta,\phi}(x)=\alpha_{N,\beta}^{-\frac{1}{\gamma_{N,\beta}}}\Big[\phi\Big(N\log\frac1{|x|}\Big)\Big]^{1-\beta},\;\;\mbox{for all}\;\; x\in B\setminus\{0\}.
\end{equation}
Then $v_{\beta,\phi}\in W^{1,N}_{0,\mathrm{rad}}(B,w_\beta)$ and
\begin{equation}\label{eq:norm}
\|v_{\beta,\phi}\|_{w_\beta}^N=(1-\beta)\int_0^{T}\phi'(t)^N\Big(\frac{t^{N-1}}{\phi(t)^N}\Big)^\beta\dd t.
\end{equation}
\end{lemma}

\begin{proof}
Let $\eta:[0,\infty)\to\R$ be a Lipschitz function such that $\eta(0)=0$ and $\eta^{\prime}\in L^{N}(0,\infty)$. Define 
\[
u_\eta(x)=\alpha_{N,\beta}^{-\frac{1}{\gamma_{N,\beta}}}\,\eta\Big(\Big(N\log\frac1{|x|}\Big)^{1-\beta}\Big),\;\;\mbox{for}\;\; x\in B\setminus\{0\}.
\]
By setting $s=(N\log\frac1r)^{1-\beta}$, with $r=|x|$, and using
the definitions of $\alpha_{N,\beta}$ and $\gamma_{N,\beta}$, we obtain
\begin{equation}\label{eq:isometry}
\int_B|\nabla u_\eta|^Nw_\beta\dd x=\int_0^\infty|\eta'(s)|^N\dd s.
\end{equation}
Note that  $\eta(s)=\big[\phi\big(s^{\frac1{1-\beta}}\big)\big]^{1-\beta}$ is a Lipschitz function on $[0, \infty)$ and that $v_{\beta,\phi}=u_\eta$. In fact, since $\phi(t)=(y_1/t_1)t$ on $[0,t_1]$, we have $\eta(s)=(y_1/t_1)^{1-\beta}s$ on $[0,t_1^{1-\beta}]$; moreover $\eta$ is constant on $[T^{1-\beta},\infty)$, and it is Lipschitz on $[t_1^{1-\beta},T^{1-\beta}]$, because $\phi$ is Lipschitz and $\phi\ge y_1>0$ there. Thus $\eta$ is Lipschitz, $\eta(0)=0$ and $\eta'\in L^N(0,\infty)$. Choose $h_k\in C_0^\infty(0,\infty)$ with $h_k\to\eta'$ in $L^N(0,\infty)$ and let $\eta_k(s)=\int_0^sh_k(\tau)\dd\tau$. Then $\eta_k$ vanishes near $s=0$ and is constant for large $s$, so $u_{\eta_k}$ vanishes near $\partial B$, is constant near the origin, and belongs to $C^\infty_{0,\mathrm{rad}}(B)$. By \eqref{eq:isometry},
\[
\|u_{\eta_k}-v_{\beta,\phi}\|_{w_\beta}^N=\|h_k-\eta'\|_{L^N(0,\infty)}^N\to0,
\]
so $v_{\beta,\phi}\in W^{1,N}_{0,\mathrm{rad}}(B,w_\beta)$. Finally, the substitution $s=t^{1-\beta}$ gives 
\[
\|v_{\beta,\phi}\|_{w_\beta}^N=\int_0^\infty|\eta'(s)|^N\dd s=(1-\beta)\int_0^\infty\phi'(t)^N\Big(\frac{t^{N-1}}{\phi(t)^N}\Big)^{\beta}\dd t,
\]
which is \eqref{eq:norm}, since $\phi'=0$ on $(T,\infty)$.
\end{proof}

\begin{lemma}\label{lem:normalized}
Let $\phi$ and $v_{\beta,\phi}$ be as in Lemma~\ref{lem:test}, and set
\begin{equation}\label{eq:theta}
\theta_\beta(\phi)=\|v_{\beta,\phi}\|_{w_\beta}^{\gamma_{N,\beta}}\;\; \mbox{and}\;\; u_{\beta,\phi}=\frac{v_{\beta,\phi}}{\|v_{\beta,\phi}\|_{w_\beta}}.
\end{equation}
Then $\|u_{\beta,\phi}\|_{w_\beta}=1$ and
\begin{equation}\label{eq:lower-beta}
\MT_N(\alpha_{N,\beta},\beta)\ge J_{N,\beta}(u_{\beta,\phi})=|B|\int_0^\infty e^{\frac{\phi(t)^{\frac{N}{N-1}}}{\theta_\beta(\phi)}-t}\dd t.
\end{equation}
\end{lemma}

\begin{proof}
Since $\phi(t)=(y_1/t_1)t$ on $[0,t_1]$, \eqref{eq:norm} gives
\begin{equation}\label{eq:first-piece}
\|v_{\beta,\phi}\|_{w_\beta}^N\ge(1-\beta)\Big(\frac{y_1}{t_1}\Big)^{N(1-\beta)}\int_0^{t_1}t^{-\beta}\dd t=\Big(\frac{y_1^N}{t_1^{N-1}}\Big)^{1-\beta}>0,
\end{equation}
so $u_{\beta,\phi}$ is well defined and $\|u_{\beta,\phi}\|_{w_\beta}=1$. By \eqref{eq:test}, \eqref{eq:theta} and $(1-\beta)\gamma_{N,\beta}=\frac N{N-1}$, we have
\[
\alpha_{N,\beta}|u_{\beta,\phi}(r)|^{\gamma_{N,\beta}}=\frac{\phi(t)^{(1-\beta)\gamma_{N,\beta}}}{\|v_{\beta,\phi}\|_{w_\beta}^{\gamma_{N,\beta}}}=\frac{\phi(t)^{\frac{N}{N-1}}}{\theta_\beta(\phi)},\;\;\mbox{with}\;\; t=N\log\frac1r, \;\;\mbox{for all}\;\; 0<r<1.
\]
Hence,
\[
J_{N,\beta}(u_{\beta,\phi})=\omega_{N-1}\int_0^1e^{\alpha_{N,\beta}|u_{\beta,\phi}(r)|^{\gamma_{N,\beta}}}r^{N-1}\dd r=|B|\int_0^\infty e^{\frac{\phi(t)^{\frac{N}{N-1}}}{\theta_\beta(\phi)}-t}\dd t.
\]
Since $u_{\beta,\phi}$ is admissible in \eqref{eq:sup}, we get 
\eqref{eq:lower-beta}.
\end{proof}

\begin{lemma}\label{lem:limit}
Let $\phi$ be a profile with nodes $t_1$ and $T$. Then
\begin{equation}\label{eq:theta-limit}
\theta_1(\phi):=\lim_{\beta\to1^-}\theta_\beta(\phi)=\Big[\frac{y_1^N}{t_1^{N-1}}\exp\Big(\int_{t_1}^{T}\frac{t^{N-1}\phi'(t)^N}{\phi(t)^N}\dd t\Big)\Big]^{\frac1{N-1}}.
\end{equation}
\end{lemma}

\begin{proof}
Let $E_1=y_1^N/t_1^{N-1}$. From \eqref{eq:first-piece} and \eqref{eq:norm}, we can write
\[
\|v_{\beta,\phi}\|_{w_\beta}^N=E_1^{1-\beta}+(1-\beta)R_\beta,\;\;\mbox{where}\;\;
R_\beta=\int_{t_1}^{T}\phi'(t)^N\Big(\frac{t^{N-1}}{\phi(t)^N}\Big)^\beta\dd t.
\]
Since $\phi\ge y_1>0$ on $[t_1,T]$ and $\phi'$ is bounded, the dominated convergence theorem gives $R_\beta\to R_1:=\int_{t_1}^{T}t^{N-1}\phi'^N\phi^{-N}\dd t$ as $\beta\to1^-$. Since $E_1^{1-\beta}=1+(1-\beta)\log E_1+O((1-\beta)^2)$, we obtain
\[
\log\theta_\beta(\phi)=\frac{\log\|v_{\beta,\phi}\|_{w_\beta}^N}{(N-1)(1-\beta)}=\frac{\log\big(1+(1-\beta)(\log E_1+R_1)+o(1-\beta)\big)}{(N-1)(1-\beta)}\to\frac{\log E_1+R_1}{N-1},
\]
which is \eqref{eq:theta-limit}.
\end{proof}

\begin{proposition}\label{prop:uniform} Let $\phi$ be a profile with nodes $t_1$ and $T$. Then,  
\begin{equation}\label{eq:FN}
\MT_N(\alpha_{N,\beta},\beta)\ge|B|\int_0^\infty e^{\frac{\phi(t)^{\frac{N}{N-1}}}{\theta_1(\phi)}-t}\dd t,
\end{equation}
 for every $\beta\in[0,1)$.
\end{proposition}

\begin{proof}
Fix $\beta\in[0,1)$ and let $\beta'\in(\beta,1)$. By \eqref{eq:mono} and Lemma~\ref{lem:normalized},
\[
\MT_N(\alpha_{N,\beta},\beta)\ge\MT_N(\alpha_{N,\beta'},\beta')\ge|B|\int_0^\infty e^{\frac{\phi(t)^{\frac{N}{N-1}}}{\theta_{\beta'}(\phi)}-t}\dd t.
\]
By Lemma~\ref{lem:limit}, $\theta_{\beta'}(\phi)\to\theta_1(\phi)>0$ as $\beta'\to1^-$, and the integrands are bounded by $\exp({\phi(T)^{N/(N-1)}/\theta_{\beta'}(\phi)-t})\le Ce^{-t}$, for $\beta'$ close to $1$. The dominated convergence theorem gives the claim.
\end{proof}

\section{The test profile}\label{sec:test-profile}
\begin{figure}[htbp]
\centering
\begin{tikzpicture}[x=0.55cm,y=0.3cm,>=stealth,every node/.style={font=\small}]
  \draw[->] (0,0) -- (17.5,0) node[right] {$t$};
  \draw[->] (0,0) -- (0,14.8) node[above] {$\phi_3(t)$};
  \node[below left] at (0,0) {$0$};
  \draw (10,0.25) -- (10,-0.25) node[below] {$10$};
  \draw (14,0.25) -- (14,-0.25) node[below] {$14$};
  \draw (0.12,10) -- (-0.12,10) node[left,font=\scriptsize] {$10$};
  \draw (0.12,12.5146) -- (-0.12,12.5146) node[left,font=\scriptsize] {$10\big(\frac75\big)^{2/3}$};
  \draw[densely dashed,gray!60] (10,0) -- (10,10) -- (0,10);
  \draw[densely dashed,gray!60] (14,0) -- (14,12.5146) -- (0,12.5146);
  \draw[very thick,blue!60!black] (0,0) -- (10,10);
  \draw[very thick,blue!60!black,domain=10:14,samples=40] plot (\x,{10*pow(\x/10,2/3)});
  \draw[very thick,blue!60!black] (14,12.5146) -- (17.2,12.5146);
  \fill[blue!60!black] (10,10) circle (2pt);
  \fill[blue!60!black] (14,12.5146) circle (2pt);
\end{tikzpicture}
\caption{The profile $\phi_3$ of \eqref{eq:phiN}: it is linear on $[0,10]$, equal to $10\,(t/10)^{2/3}$ on $[10,14]$ and constant on $[14,\infty)$.}
\label{fig:phi3}
\end{figure}
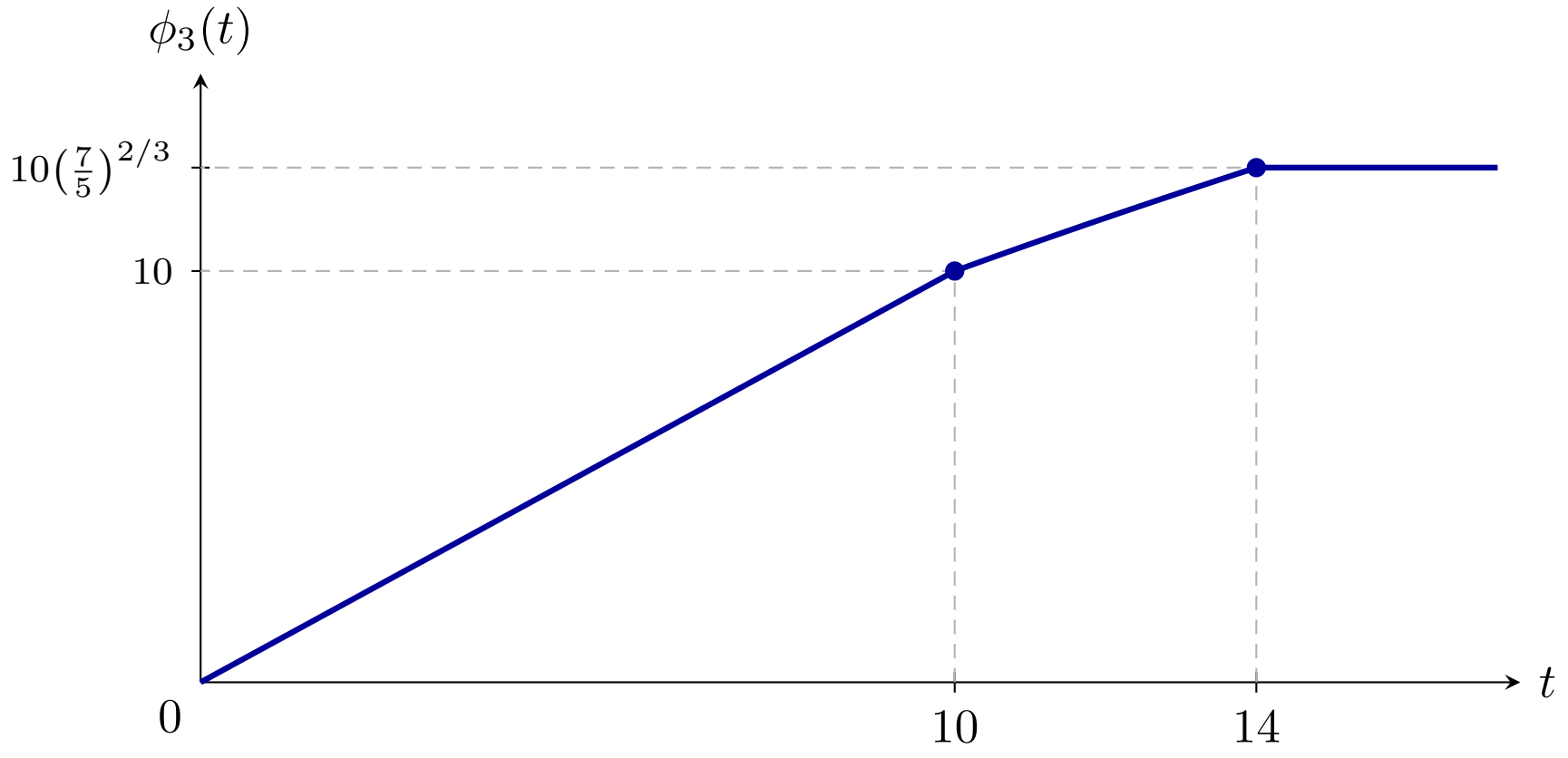

Throughout this section $N\ge3$ is fixed. Let $t_1=2(2N-1)$ and $T=\frac75\,t_1$, and let $\phi_N$ be the profile with nodes $t_1$ and $T$ given by
\begin{equation}\label{eq:phiN}
\phi_N(t)=\begin{cases}
t,&0\le t\le t_1,\\[2pt]
t_1\Big(\dfrac t{t_1}\Big)^{\frac{N-1}N},&t_1\le t\le T,\\[10pt]
t_1\Big(\dfrac75\Big)^{\frac{N-1}N},&t\ge T.
\end{cases}
\end{equation}
The function $\phi_N$ is concave and $1$-Lipschitz, and $y_1=t_1$ in Definition~\ref{def:profile}; see Figure~\ref{fig:phi3}. Like the test functions of Carleson and Chang \cite[proof of Step~II]{CarlesonChang86} (see also \cite[Section~2]{HudsonLeckband02}), $\phi_N$ is linear on an initial interval and is then given by a power. The exponent $\frac{N-1}N$ is chosen so that $t\phi_N'(t)/\phi_N(t)$ is constant on $(t_1,T)$. As the next lemma shows, this makes the integrand in \eqref{eq:FN} with $\phi=\phi_N$ an explicit exponential on $[t_1,\infty)$. 

\begin{lemma}\label{lem:explicit}
Let $\phi_N$ be given by \eqref{eq:phiN} and let
\begin{equation}\label{eq:delta}
\delta=\delta_N:=1-\Big(\frac57\Big)^{\frac1{N-1}\big(1-\frac1N\big)^N}.
\end{equation}
Then
\begin{equation}\label{eq:FN-explicit}
\int_0^\infty e^{\frac{\phi_N(t)^{\frac{N}{N-1}}}{\theta_1(\phi_N)}-t}\dd t=\int_0^{t_1}e^{(1-\delta)\,t\,(\frac{t}{t_1})^{\frac1{N-1}}-t}\dd t+\int_{t_1}^{T}e^{-\delta t}\dd t+e^{-\delta T}.
\end{equation}
\end{lemma}

\begin{proof}
On $(t_1,T)$ we have
\[
\frac{t^{N-1}\phi_N'(t)^N}{\phi_N(t)^N}=\Big(1-\frac1N\Big)^N\frac1t.
\]
Since $y_1=t_1$, \eqref{eq:theta-limit} gives
\[
\theta_1(\phi_N)^{N-1}=t_1\Big(\frac75\Big)^{(1-\frac1N)^N},
\;\;\text{that is,}\;\;
\frac1{\theta_1(\phi_N)}=(1-\delta)\,t_1^{-\frac1{N-1}}.
\]
Hence, we obtain
\[
 \frac{\phi_N(t)^{\frac{N}{N-1}}}{\theta_1(\phi_N)}=\left\{ \begin{aligned}
&(1-\delta)\,t\,\Big(\frac{t}{t_1}\Big)^{\frac1{N-1}},\;\;&\mbox{if}&\;\; 0\le t\le t_1,\\
&(1-\delta)t,\;\;&\mbox{if}&\;\; t_1\le t\le  T,\\
&(1-\delta)\,T,\;\;&\mbox{if}&\;\; t\ge  T.
  \end{aligned}\right.
\]
 Since $\int_T^\infty e^{(1-\delta)T-t}\dd t=e^{-\delta T}$, this gives \eqref{eq:FN-explicit}.
\end{proof}
The only property of $\delta$ that we need is the following upper bound.
\begin{lemma}\label{lem:delta}
The number $\delta$ in \eqref{eq:delta} satisfies $4(2N-1)\delta<1$.
\end{lemma}

\begin{proof}
Recalling that  $1-e^{-y}\le y$, for all $y\in\mathbb R$, we have
\[
\delta=1-\Big(\frac57\Big)^{\frac1{N-1}\left(1-\frac1N\right)^N}=1-e^{-\frac{1}{N-1}\left(1-\frac1N\right)^N\log\big(\frac75\big)} \le\frac1{N-1}\big(1-\frac1N\big)^N\log\Big(\frac75\Big).
\]
Hence it suffices to show that
\begin{equation}\label{eq:kappa-bound}
8\log\Big(\frac75\Big)<\frac{2(N-1)}{2N-1}\Big(\frac N{N-1}\Big)^N.
\end{equation}
Now, we claim that
\begin{equation}\label{claim-elem}
 \log\frac{1+x}{1-x}\ge2x+\frac{2x^3}3\;\;\text{and}\;\;
\log(1-x)\ge-x-\frac{x^2}2-\frac{x^3}3-\frac{x^4}{4(1-x)}, 
\end{equation}
for all $x\in [0, 1)$. Indeed, by setting $f,g:[0,1)\to\mathbb{R}$ given by
\[
f(x) = \log(1+x) - \log(1-x) - 2x - \frac{2x^3}{3}
\]
and
\[
g(x) = \log(1-x) + x + \frac{x^2}{2} + \frac{x^3}{3} + \frac{x^4}{4(1-x)},
\]
we can see that
\[
f'(x)= \frac{2x^4}{1-x^2}\;\;\mbox{and}\;\; g'(x) = \frac{x^4}{4(1-x)^2}.
\]
Since $f(0)=g(0)=0$, $f^{\prime}(x)\ge 0$ and $g^{\prime}(x)\ge 0$ on $[0, 1)$, we obtain \eqref{claim-elem}. Let $x=\frac1{2N-1}\in(0,\frac15]$, for $N\ge 3$. Since
\[
\frac N{N-1}=\frac{1+x}{1-x},\;\; N=\frac{1+x}{2x}\;\;\text{and}\;\;\frac{2(N-1)}{2N-1}=1-x,
\]
the claim \eqref{claim-elem} yields
\begin{align*}
  \log\Big[\frac{2(N-1)}{2N-1}\Big(\frac N{N-1}\Big)^N\Big]
&=\log(1-x)+\frac{1+x}{2x}\log\frac{1+x}{1-x}\\
&\ge1-\frac{x^2}6-\frac{x^4}{4(1-x)}\\
&\ge1-\frac1{150}-\frac1{2000}=1-\frac{43}{6000}.
\end{align*}
Since $e>2.718$ and $e^{-y}\ge1-y$ for $y\ge 0$, we get
\begin{equation}\label{baixo}
  \frac{2(N-1)}{2N-1}\Big(\frac N{N-1}\Big)^N\ge e\Big(1-\frac{43}{6000}\Big)>2.698.
\end{equation}
On the other hand, $\frac75=\frac{1+1/6}{1-1/6}$, so from the series expansion of $\log(\frac{1+1/6}{1-1/6})$ we get
\begin{equation}\label{cima}
  \log\Big(\frac75\Big)=2\sum_{i\ge0}\frac{6^{-2i-1}}{2i+1}\le2\Big(\frac16+\frac{6^{-3}}3\sum_{i\ge0}6^{-2i}\Big)=2\Big(\frac16+\frac1{630}\Big)=\frac{106}{315}
\end{equation}
and $8\log\big(\frac75\big)\le\frac{848}{315}<2.693$. Thus, from \eqref{baixo} and \eqref{cima}, we obtain \eqref{eq:kappa-bound}.
\end{proof}

\section{Proofs of the main results}\label{sec:proofs}
\begin{proof}[Proof of Theorem~\ref{thm:main}]
Let $N\ge3$ and $\beta\in[0,1)$, and let $t_1=2(2N-1)$, $T=(7/5)t_1$, $\phi_N$ and $\delta$ be as in \eqref{eq:phiN} and \eqref{eq:delta}. By Proposition~\ref{prop:uniform} and Lemma~\ref{lem:explicit}, we get
\begin{equation}\label{eq:FN-expl2}
  \frac{\MT_N(\alpha_{N,\beta},\beta)}{|B|}\ge \int_0^{t_1}e^{(1-\delta)\,t\,(\frac{t}{t_1})^{\frac1{N-1}}-t}\dd t+\int_{t_1}^{T}e^{-\delta t}\dd t+e^{-\delta T}.
\end{equation}
We estimate the three terms on the right-hand side of \eqref{eq:FN-expl2}. First, recalling that
$t_1=2(2N-1)$,
we obtain
\[
\frac1{t_1}\int_0^{t_1}\Big[(1-\delta)\,t\Big(\frac t{t_1}\Big)^{\frac1{N-1}}-t\Big]\dd t=(1-\delta)\,\frac{N-1}{2N-1}\,t_1-\frac{t_1}2=-1-2(N-1)\delta.
\]
By Lemma~\ref{lem:delta}, we can write
\[
2(N-1)\delta<\frac{N-1}{2(2N-1)}=\frac14-\frac1{4(2N-1)}.
\]
Hence, by Jensen's inequality and $e^y\ge1+y$, for $y\in\mathbb{R}$,
\begin{equation}\label{eq:first-term}
  \begin{aligned}
  \int_0^{t_1} e^{(1-\delta)t\left(\frac{t}{t_1}\right)^{\frac{1}{N-1}} - t} \dd t &\ge t_1 e^{\frac{1}{t_1} \int_0^{t_1} \big[(1-\delta)t\left(\frac{t}{t_1}\right)^{\frac{1}{N-1}} - t \big] \dd t}\\
  &=2(2N-1)e^{-1-2(N-1)\delta}\\
  &\ge 2(2N-1)e^{-\frac{5}{4}+\frac{1}{4(2N-1)}}\\
  &\ge 2(2N-1)\,e^{-\frac54}\Big(1+\frac1{4(2N-1)}\Big)\\
  &=2e^{-\frac54}\Big(2N-\frac34\Big).
  \end{aligned}
\end{equation}
In addition, from Lemma~\ref{lem:delta} we have $-\delta>-\frac{1}{4(2N-1)}$. Hence
\begin{equation}\label{eq:second-term}
  \begin{aligned}
  \int_{t_1}^Te^{-\delta t}\dd t\ge\int_{t_1}^Te^{-\frac{t}{4(2N-1)}}\dd t=4(2N-1)\big(e^{-\frac12}-e^{-\frac7{10}}\big)
  \end{aligned}
\end{equation}
and 
\begin{equation}\label{eq:third-term}
  \begin{aligned}
  e^{-\delta T}\ge e^{-\frac{T}{4(2N-1)}}=e^{-\frac{7}{10}}.
  \end{aligned}
\end{equation}
Combining \eqref{eq:FN-expl2} with \eqref{eq:first-term}, \eqref{eq:second-term} and \eqref{eq:third-term}, we get
\begin{equation}\label{MT-cotab}
    \MT_N(\alpha_{N,\beta},\beta)\ge|B|\Big[2e^{-\frac54}\Big(2N-\frac34\Big)+4(2N-1)\,e^{-\frac12}-(8N-5)\,e^{-\frac7{10}}\Big].
\end{equation}
Now, we observe that 
\begin{equation}\label{e>estimate}
e^{-\frac{7}{10}} < \frac{497}{1000},\;\;\; e^{-\frac12} > \frac{303}{500}\;\;\mbox{and}\;\; e^{-\frac54} > \frac{143}{500}.
\end{equation}
Indeed, applying the second inequality in \eqref{claim-elem} with $x = \frac{18}{61}$ gives
\[
\log\left(\frac{43}{61}\right) \ge -\frac{18}{61} - \frac{162}{3721} - \frac{1944}{226981} - \frac{26244}{9760183} = -\frac{3414816}{9760183} > -\frac{7}{20},
\]
which implies $e^{-7/20} < \frac{43}{61}$, and hence
\[
e^{-\frac{7}{10}} = \left(e^{-\frac{7}{20}}\right)^2 < \left(\frac{43}{61}\right)^2 = \frac{1849}{3721} < \frac{497}{1000}.
\]
Moreover, applying the first inequality in \eqref{claim-elem} with $x = \frac{25}{201}$ gives
\[
\log\left(\frac{113}{88}\right) \ge \frac{50}{201} + \frac{2}{3}\left(\frac{25}{201}\right)^3 = \frac{6091400}{24361803} > \frac{1}{4},
\]
which implies $e^{-1/4} > \frac{88}{113}$. Hence, using also that $e<2.72$,
\[
e^{-\frac{1}{2}} = \left(e^{-\frac{1}{4}}\right)^2 > \left(\frac{88}{113}\right)^2 = \frac{7744}{12769} > \frac{303}{500}
\quad\text{and}\quad
e^{-\frac{5}{4}} = e^{-1}e^{-\frac{1}{4}} > \frac{25}{68}\cdot\frac{88}{113} = \frac{550}{1921} > \frac{143}{500}.
\]
Consequently, the expression in brackets of \eqref{MT-cotab} is larger than
\begin{equation}\label{cota-menor}
    \begin{aligned}
       \frac{143}{250}\Big(2N-\frac34\Big)+\frac{606}{250}(2N-1)-\frac{497}{1000}(8N-5) & = \frac{143+606-497}{125}\,N - \frac{429+2424-2485}{1000}\\
       &=\frac{252}{125}N - \frac{46}{125} > 2N-\frac25.
    \end{aligned}
\end{equation}
Thus, from \eqref{MT-cotab} and \eqref{cota-menor}
\[
\MT_N(\alpha_{N,\beta},\beta)>|B|\Big(2N-\frac{2}{5}\Big).
\]
It remains to compare $2N-\frac25$ with the concentration level. Precisely, we need to prove that 
\begin{equation}\label{lower bound for CL}
   1+e^{1+\frac12+\cdots+\frac1{N-1}}<2N-\frac25,\;\;\mbox{for all}
\;\; N\ge 3. 
\end{equation}
\begin{remark}
This is the step of the proof in which the hypothesis $N\ge3$ is essential. In fact, for $N=2$, the left-hand side of \eqref{lower bound for CL} equals $1+e>\frac{18}5=2\cdot2-\frac25$.
\end{remark}
 In order to obtain \eqref{lower bound for CL}, we first suppose that $N=3$. Applying the first inequality in \eqref{claim-elem} with $x = \frac{4}{11} \in [0,1)$, we obtain
\[
\log\left(\frac{15}{7}\right) \ge \frac{8}{11} + \frac{2}{3}\left(\frac{4}{11}\right)^3 = \frac{3032}{3993} > \frac{3}{4}.
\]
Consequently, $2\log\left(\frac{15}{7}\right)>\frac{3}{2}=1+\frac{1}{2}$ and we have
\[
1+e^{1+\frac{1}{2}}<1+\left(\frac{15}{7}\right)^2 = 1+\frac{225}{49} < 1+\frac{23}{5} = 2\cdot 3-\frac25.
\]
Let $N\ge4$.  From \eqref{claim-elem}, we  have 
\[
\int_{k-\frac{1}{2}}^{k+\frac12}\frac{\dd x}x=\log\frac{1+\frac{1}{2k}}{1-\frac{1}{2k}}\ge \frac{1}{k},\;\;\mbox{for}\;\; k\ge 1.
\]
Therefore
\[
  \begin{aligned}
    \frac13+\frac{1}{4}+\cdots+\frac1{N-1}&\le \int_{3-\frac12}^{3+\frac12}\frac{\dd x}{x}+\int_{4-\frac12}^{4+\frac12}\frac{\dd x}{x}+\cdots+\int_{N-1-\frac12}^{N-1+\frac12}\frac{\dd x}{x}\\
    &=\int_{\frac{5}{2}}^{\frac{2N-1}{2}}\frac{\dd x}{x}=\log\frac{2N-1}{5}.
  \end{aligned}
\]
Thus, since $e^{1+\frac12}<\frac{23}{5}$ by the case $N=3$, we get
\[
1+e^{1+\frac12+\cdots+\frac1{N-1}}\le1+e^{1+\frac12}\,\frac{2N-1}5<1+\frac{23}{25}\,(2N-1)=2N-\frac25-\frac{4(N-3)}{25}\le2N-\frac25,
\]
which is \eqref{lower bound for CL}. Hence
\[
\MT_N(\alpha_{N,\beta},\beta)>\big(2N-\frac25\big)|B|>|B|\big(1+e^{1+\frac12+\cdots+\frac1{N-1}}\big).
\]
\end{proof}
\begin{remark}\label{rem:profile-lower-bound}
The proof of Theorem~\ref{thm:main} gives more than its statement. Indeed, by Proposition~\ref{prop:uniform} and Lemma~\ref{lem:explicit} (see \eqref{eq:FN-expl2}), for every $\beta\in[0,1)$ we have the sharper lower bound
\begin{equation}\label{Precise:FN-expl2}
  \frac{\MT_N(\alpha_{N,\beta},\beta)}{|B|}\ge F_N:=\int_0^{t_1}e^{(1-\delta_N)\,t\,(\frac{t}{t_1})^{\frac1{N-1}}-t}\dd t+\int_{t_1}^{T}e^{-\delta_N t}\dd t+e^{-\delta_N T},
\end{equation}
where
\[
t_1=2(2N-1),\quad T=\frac75\,t_1\quad\text{and}\quad \delta_N=1-\Big(\frac57\Big)^{\frac1{N-1}\left(1-\frac1N\right)^N}.
\]
The estimates in the proof of Theorem~\ref{thm:main} show that $F_N>2N-\frac25$. Thus, the constant $2N-\frac25$ in Theorem~\ref{thm:main} is a simple explicit lower bound for $F_N$, and not the full estimate given by the profile $\phi_N$. The explicit expression for $F_N$, which does not depend on $\beta$, may also serve as a benchmark when other profiles are considered or when improved quantitative bounds are sought.
\end{remark}
\begin{proof}[Proof of Theorem~\ref{thm:attainment}]
  Fix $N\ge3$ and $\beta\in[0,1)$. From \eqref{eq:concentration} and Theorem~\ref{thm:main},
\begin{equation}\label{eq:gap}
\MT_N(\alpha_{N,\beta},\beta)>|B|\big(1+e^{1+\frac12+\cdots+\frac1{N-1}}\big)\ge\Sigma_{N,\beta}.
\end{equation}
 Let $(u_k)$ be a maximizing sequence for $\MT_N(\alpha_{N,\beta},\beta)$. Since $t\mapsto J_{N,\beta}(tu_k)$ is nondecreasing on $[0,\infty)$, we may assume that $\|u_k\|_{w_\beta}=1$. Up to a subsequence, $u_k\weakto u_0$ weakly in $W^{1,N}_{0,\mathrm{rad}}(B,w_\beta)$. By the concentration-compactness alternative \cite[Lemma~2.1]{RoyN19}, up to a further subsequence, either $J_{N,\beta}(u_k)\to J_{N,\beta}(u_0)$, or $\int_{B\setminus B_r}|\nabla u_k|^Nw_\beta\dd x\to0$ for every $r\in(0,1)$. In the second case, since $u\mapsto\int_{B\setminus B_r}|\nabla u|^Nw_\beta\dd x$ is convex and continuous on $W^{1,N}_{0,\mathrm{rad}}(B,w_\beta)$, and hence weakly lower semicontinuous, we get $\int_{B\setminus B_r}|\nabla u_0|^Nw_\beta\dd x=0$ for every $r\in(0,1)$; letting $r\to0$, we obtain $\|u_0\|_{w_\beta}=0$, that is, $u_0\equiv0$. Thus $(u_k)$ is a normalized concentrating sequence and $\MT_N(\alpha_{N,\beta},\beta)=\lim_kJ_{N,\beta}(u_k)\le\Sigma_{N,\beta}$, which contradicts \eqref{eq:gap}. Hence $J_{N,\beta}(u_0)=\lim_kJ_{N,\beta}(u_k)=\MT_N(\alpha_{N,\beta},\beta)$. Moreover, $\|u_0\|_{w_\beta}\le\liminf_k\|u_k\|_{w_\beta}=1$, so $u_0$ is a maximizer.
\end{proof}

\noindent\textbf{Funding.} J.\,F.~de Oliveira acknowledges partial support from the National Council for Scientific and Technological Development (CNPq, Brazil) through grant 303443/2025-1. 
J. M. do \'O acknowledges partial support from CNPq
through grants 312435/2026-6, 409764/2023-0, and 443594/2023-6, CAPES MATH AMSUD grant 88887.878894/2023-00 and Para\'iba State Research Foundation (FAPESQ), grant no 3034/2021. R.\,C.~Ponciano acknowledges partial support from São Paulo Research Foundation (FAPESP) grant 2023/07697-9.

\medskip

\noindent{\bf AI Declaration.} All mathematical content of this paper, including the statements of the results and their proofs, was developed by the authors, who take full responsibility for the entire content of the manuscript. After the mathematical work was completed, AI-assisted tools were used solely for proofreading and for improving the language and readability of the text. All suggestions produced by these tools were reviewed and verified by the authors.

\medskip

\noindent\textbf{Data availability statement:} Our manuscript has no associated data or further material.

\medskip

\noindent\textbf{Ethical Approval:}  Not applicable.

\medskip

\noindent\textbf{Conflict of interest.} The authors declare no conflict of interest.


\end{document}